\documentstyle{amsart}

\newcommand{\lan}{\langle}
\newcommand{\ran}{\rangle}

\newcommand{\Card}{\operatorname{Card}}
\newcommand{\card}{\operatorname{card}}
\newcommand{\Range}{\operatorname{Range}}
\newcommand{\ORD}{\operatorname{ORD}}

\newsymbol\restriction 1316
\newsymbol\Vdash 130D

\newcommand{\vDash}{\models}
\begin{document}

\Large

\baselineskip=.30truein

\centerline{\bf Coding without Fine Structure}

\vskip30pt

\large
   \begin{center}

Sy D. Friedman\footnote{Research supported by NSF Contract \# 9205530.}
\end{center}

  \begin{center}
M.I.T.\end{center}

\vskip10pt

In this paper we prove Jensen's Coding Theorem, assuming $\sim 0^{\#}$,
via a proof that makes no use of the fine structure theory. We do need to
quote Jensen's Covering Theorem, whose proof uses fine-structural ideas,
but make no direct use of these ideas. The key to our proof is the use of
``coding delays.''

\noindent
{\bf Coding Theorem} \ (Jensen) \ Suppose $\lan M,A\ran$ is a model of
$ZFC + O^{\#}$ does not exist. Then there is an $\lan M,A\ran$-definable
class forcing $P$ such that if $G\subseteq P$ is $P$-generic over $\lan
M,A\ran:$

(a) \ $\lan M[G], A,G\ran \vDash ZFC.$

(b) \ $M[G]\vDash V=L[R], R\subseteq \omega $ and $\lan M[G],A,G\ran \vDash
A,G$ are definable from the parameter $R$.

In the above statement when we say ``$\lan M,A\ran\vDash ZFC$'' we mean
that $M\vDash ZFC$ and in addition $M$ satisfies replacement for formulas
that mention $A$ as a predicate. And ``$P$-generic over $\lan M,A\ran$''
means that all $\lan M,A\ran$-definable dense classes are met.

The consequence of $\sim O^{\#}$ that we need follows directly from the
Covering Theorem.

\vskip10pt

\noindent
{\bf Covering Theorem}\ (Jensen) \ Assume $\sim O^{\#}$.  If $X$ is an
uncountable set of ordinals then there is a constructible $Y\supseteq X,$
card $Y=\operatorname{card} X.$

\vskip10pt

\noindent
{\bf Lemma 1}\ (Jensen) \ Assume $\sim O^{\#}$. If $j:L_\alpha
\longrightarrow L_\beta $ is $\Sigma _1$-elementary, $\alpha \ge\omega _2$
and $\kappa =crit(j)$ then $\alpha <(\kappa ^+)^L.$

\vskip10pt

\noindent
{\bf Proof} \ Of course $crit(j)$ denotes the least ordinal $\kappa $ such
that $j(\kappa )\neq\kappa,$ which we assume to exist. Now let
$U=\{X\subseteq \kappa |X\in L_\alpha ,\kappa \in j(X)\}.$ If $\alpha \ge
(\kappa ^+)^L$ then $U$ is an ultrafilter on all constructible subsets of
$\kappa $ and we can form $Ult(L,U)=$ ultrapower of $L$ by $U$ (using
constructible functions to form the ultrapower).  If this is well-founded
then we get a nontrivial elementary embedding $L\longrightarrow L,$ which
gives $O^{\#}$ by a theorem of Kunen.

Now we know that $Ult(L_\alpha ,U)$ is well-founded since it embeds into
$L_\beta $ (using: $k([f])=j(f)(\kappa )).$ And by a Lowenheim-Skolem
argument, if $Ult(L,U)$ were ill-founded then so would be $Ult(L_{\kappa
^+},U),\kappa ^+=$ the real $\kappa ^+.$ 
So we may assume that $\kappa \ge \omega _2$ as otherwise $\kappa ^+\le
\omega _2\le \alpha $ and the facts above would imply that
$Ult(L,U)$ were well-founded.

Using the Covering Theorem and the fact that $\kappa \ge\omega _2$ we show
that if $\lan X_n|n\in\omega \ran$ belong to $U$ then $\bigcap\limits_{n}$
$X_n\neq \phi $ ($U$ is ``countably complete''), a fact that immediately
yields the well-foundedness of $Ult(L,U).$

Apply Covering to get $F\in L$ of cardinality $\omega _1$ such that $X_n\in
F$ for each $n.$ As $\kappa \ge\omega _2,F$ has $L$-cardinality $<\kappa $
and also we may assume that $F$ is a subset of $P(\kappa )\cap L.$ So $F\in
L_{(\kappa ^+)^L}\subseteq L_\alpha$ and there is a bijection $h:F\longrightarrow\gamma, \gamma < \kappa, h\in L_\alpha.$ Let $F^*=\{X\in F|\kappa \in \j(X) \};$
then $F^*\in L_\alpha$ since $h[F^*] = \{j(h)(Y)|Y \in j(F), \kappa \in Y \}$
belongs to $L_\beta$ and hence to $L_\kappa \subseteq L_\alpha.$ So $\cap F^*\neq\phi $ since $j(\cap
F^*)=\cap j[F^*]$ contains $\kappa $ and $j$ is $\Sigma _1$-elementary. As
$\{X_n|n\in\omega \}\subseteq F^*$ we get $\bigcap\limits_{n}X_n\neq\phi ,$
as desired. $\dashv$

Next we show that to prove the Coding Theorem we may assume that the GCH
holds in $M,$ and that instead of coding into a real, it is enough to code
into a ``reshaped'' subset of $\omega _1$.

\vskip10pt

\noindent
{\bf Lemma 2}\ (Folklore) \ If $\lan M,A\ran$ is a model of
ZFC then there is an $\lan M,A\ran$-definable forcing $P^*$ such that if
$G^*$ is $P^*$-generic over $\lan M,A\ran$ then for some $B\subseteq
\operatorname{ORD} (M),$ $B$ is definable over $\lan M[G^*],A,G^*\ran$ and
this model satisfies $ZFC+GCH+V=L[B]+A,G^*$ are definable relative to $B.$
And if $M$ satisfies $\sim O^{\#}$ then so does $M[G^*].$

\vskip10pt

\noindent
{\bf Proof} \ First, by forcing with conditions $p:\alpha \longrightarrow
2,\alpha \in ORD,$ ordered by $p\le q$ iff $p$ extends $q$ we can obtain
$B$ as above, except for the GCH. This is beause if $G^*_0$ is generic for
this forcing and $B_0=\{\beta |p(\beta )=1$ for some $p\in G^*_0\}$ then
$M[G^*_0]\vDash V=L[B_0]$ and using $B_0$ we can identify $A$ with a
class of ordinals $B_1;$ let $B=$ the join of $B_0,B_1.$

Second, we force over $\lan L[B],B\ran$ to obtain the GCH. As usual,
$\sqsupset_\alpha $
is defined (in $L[B])$ by: $\sqsupset_0=\omega ,\sqsupset_{\alpha
+1}=2^{\sqsupset_\alpha }$ 
and $\sqsupset_\lambda =\cup\{\sqsupset_\alpha |\alpha <\lambda \}$ for limit
$\lambda.$ 
For any $\alpha$ $ P(\alpha )$ is the forcing whose conditions are $p: \ \beta
\longrightarrow 2^{\sqsupset_{\alpha }},$ $\beta <\sqsupset^+_\alpha,$ ordered
by $p\le q$ iff $p$ extends $q.$ We take $P$ to be the ``Easton product'' of
the $P(\alpha )$'s: a condition in $P$ is $p: \ \alpha (p)\longrightarrow
L[B],$ $p\in L[B]$ such that $p(\alpha )\in P(\alpha )$ for each
$\alpha <\alpha (p)$ and such that $\{\beta <\alpha |p(\beta )\neq \phi \}$ is
bounded in $\alpha $ for inaccessible $\alpha \le\alpha (p).$ For any
$\alpha $ $P$ factors as
$P(>\alpha )\times P(\le \alpha )$ where $P(>\alpha )$ is
$\sqsupset^+_{\alpha +1}$-closed and $P(\le\alpha )$ has cardinality
$\le\sqsupset_{\alpha +1}.$ It follows that ZFC is preserved, the infinite
successor cardinals of the generic extension are the $\sqsupset^+_\alpha $
of $L[B]$ and that the GCH holds in the generic extension. And if $L[B]$
satisfies $\sim O^{\#}$ then so does the $P$-generic extension, since for
singular strong limit cardinals $\kappa $ of $L[B],$ $\kappa ^+$ of
$L[B]=\kappa ^+$ of $L$ and $\kappa ^+$ of $L[B]=\kappa ^+$ of the
$P$-generic extension.

Let $P^*$ be the product of the two forcings described above. $\dashv$

\vskip10pt

\noindent
{\bf Definition} \ $b\subseteq \omega_1 $ is {\it reshaped} if $\xi <\omega
_1\longrightarrow \xi $ is countable in $L[b\cap\xi ].$

\vskip10pt

\noindent
{\bf Lemma 3} \ (Jensen-Solovay [68]) \  Suppose $M\vDash ZFC+V=L[b]$ where
$b$ is a reshaped subset of $\omega _1.$ Then there is a CCC forcing $P$
such that if $G$ is $P$-generic over $M$ then $M[G]\vDash V=L[R]$ where
$R\subseteq \omega .$

\vskip10pt

\noindent
{\bf Proof} \ Using the fact that $b$ is reshaped we may choose $\lan
R'_\xi |\xi <\omega _1\ran$ so that for each $\xi <\omega _1,$ $R'_\xi $ is
the least real in $L[b\cap\xi ]$ distinct from each $R'_{\xi '},$ $\xi
'<\xi.$  Let $R_\xi =\{n<\omega |n$ codes a finite initial segment of the
characteristic function of $R'_\xi \}.$ Then $\xi _0\neq\xi
_1\longrightarrow R_{\xi _0}\cap R_{\xi _1}$ is finite.

A condition in $P$ is $p=(s(p),s^*(p))$ where $s(p)$ is a finite subset of
$\omega $ and $s^*(p)$ is a finite subset of $b.$ Extension is defined by: \
$p\le q$ iff $s(p)$ end extends $s(q),$ $s^*(p)\supseteq s^*(q)$ and $\xi
\in s^*(q)\longrightarrow s(p)-s(q)$ is disjoint from $R_\xi .$ This is ccc
and if $G$ is $P$-generic, $R=\cup\{s(p)|p\in G\}$ then $\xi \in b$ iff $R\cap
R_\xi $ is finite. So inductively we can recover $b\cap\xi ,R_\xi $ in
$L[R].$ And $p\in G$ iff $s(p)$ is an initial segment of $R,$ $\xi \in
s^*(p)\longrightarrow R_\xi \cap R\subseteq s(p).$ So $M[G]\vDash
V=L[b,G]=L[R].$ $\dashv$

\vskip10pt

Thus the Coding Theorem with $\sim O^{\#}$ reduces to:

\vskip10pt

\noindent
{\bf Theorem 4} \ Suppose that $A\subseteq \operatorname{ORD}$ and
$\lan L[A],A\ran$ is a model of ZFC $+$ GCH $+ \sim O^{\#}.$ Then there is an
$\lan L[A],A\ran$-definable class forcing $P$ such that if $G$ is
$P$-generic over $\lan L[A],A\ran:$

(a) \ $\lan L[A,G], A,G\ran$ is a model of ZFC.

(b) \ $L[A], L[A,G]$ have the same cofinalities.

(c) \ $L[A,G]=L[X]$ where $X$ is a reshaped subset of $\omega _1$ and $A,G$
are definable over $L[X]$ with parameter $X.$

It is useful to make the following harmless assumption about $A:$ if
$H_\alpha ,\alpha $ an infinite $L[A]$-cardinal, denotes $\{X\in L[A]|$
transitive closure $(X)$ has $L[A]$-cardinality $<\alpha \}$ then $H_\alpha
=L_\alpha [A].$ This is easily arranged using the GCH in $L[A].$

\vskip10pt

\noindent
{\bf Definition of the Forcing P}

Let $\Card = $ all infinite cardinals, $\Card^+=\{\alpha ^+|\alpha \in
\Card\}$ and $\Card'=$ all uncountable limit cardinals. Of course these
definitions are made in $V=L[A].$

\vskip10pt

\noindent
{\bf Definition} \ (Strings) \ Let $\alpha \in\Card.$ $S_\alpha $ consists
fo all $s: \ [\alpha ,|s|)\longrightarrow 2,$ $\alpha \le |s|<\alpha ^+$ such
that $|s|$ is a multiple of $\alpha $  and for all $\eta \le |s|,$
$L_\delta [A\cap\alpha ,s\restriction\eta ]\vDash\card(\eta )\le\alpha $
for some $\delta <(\eta ^+)^L\cup \omega _2.$

Thus for $\alpha \ge \omega _2$ we insist that $s$ is ``quickly reshaped'' in
that $\eta \le|s|$ is collapsed relative to $A\cap \alpha ,$
$s\restriction\eta $ before $(\eta ^+)^L.$ This will enable us to establish
cofinality-preservation, using Lemma 1. Note that we allow $|s|=\alpha ,$
in which case $s=\phi _\alpha,$ the ``empty string at $\alpha.$'' Also for
$s,t\in S_\alpha $ write $s\le t$ for $s\subseteq t$ and $s<t$ for $s\le
t,$ $s\neq t.$

\vskip10pt

\noindent
{\bf Definition} \ (Coding Structures) \ For $s\in S_\alpha $ defome $\mu
^{<s},\mu ^s$ inductively by: $\mu ^{<\phi _\alpha }=\alpha ,\mu
^{<s}=\cup\{\mu ^t|t<s\}$ for $s\neq \phi _\alpha $ and $\mu ^s=$ least $\mu
>\mu ^{<s}$ such that $\mu '\mu =\mu $ for $\mu '<\mu $ and $L_\mu
[A\cap\alpha ,s]\vDash$ ``$s\in S_\alpha $''.  And $\cal{A}^s=L_{\mu
^s}[A\cap\alpha ,s],$ $\cal{A}^{<s}=\lan L_{\mu ^{<s}}[A\cap\alpha ,\hat
s],A\cap\alpha ,\hat s\ran$ where $\hat s=\{\mu ^{<t}+\delta |t<s,\delta <\alpha ,s(|t|+\delta )=1\}.$

Thus by definition there is $\delta <\mu ^s$ such that $L_\delta
[A\cap\alpha ,s]$ $\vDash\card(|s|)\le\alpha $ and $L_{\mu
^s}\vDash\card(\delta )\le|s|,$ when $\alpha \ge \omega _2.$ For $|s|=\eta
+\alpha ,\eta $ a multiple of $\alpha ,\cal{A}^{<s}$ has universe
$\cal{A}^{s\restriction \eta }$ and for $|s|$ a limit of multples of $\alpha ,\cal{A}^{<s}=\cup\{\cal{A}^{<t}|t<s\}.$

\vskip10pt

\noindent
{\bf Definition} \ (Coding Apparatus) \ For $\omega \neq \alpha \in\Card,
s\in S_\alpha ,$ $i<\alpha $ let $H^s(i)=\Sigma _1$ Hull of $i\cup\{A\cap
\alpha ,s\}$ in $\cal{A}^s$ and $f^s(i)=$ ordertype $(H^s(i)\cap \ORD).$
For $\alpha \in \Card^+,b^s=\Range(f^s\restriction B^s)$ where
$B^s=\{i<\alpha |i=H^s(i)\cap\alpha \}.$ 
Also for $\eta <|s|,$ $\eta =|t|+\delta ,$ $\delta <\alpha ,t<s$ we define
$b^{s\restriction\eta }=\{\gamma +\delta |\gamma \in b^t\}.$

\vskip10pt

\noindent
{\bf Definition} \ (A Partition of the Ordinals) \ Let $B,C,D,E$ denote
the classes of ordinals congruent to $0,1,2,3$ mod $4,$ respectively.
Also for any ordinal $\alpha $ and $X=B,C,D$ or $E$ we write $\alpha ^X$
for the $\alpha ^{\text{th}}$ element of $X.$

\vskip10pt

\noindent
{\bf Definition} \ (The Successor Coding) \ Suppose $\alpha \in \Card, s\in
S_{\alpha ^+}.$ A condition in $R^s$ is a pair $(t,t^*)$ where $t\in
S_\alpha ,$ $t^*\subseteq \{b^{s\restriction\eta }|\alpha \le\eta <|s|\},$
$\card(t^*)\le\alpha.$ Extension of conditions is defined by:
$(t_0,t^*_0)\le(t_1,t^*_1)$ iff $t_1\le t_0,$ $t^*_1\subseteq t^*_0$ and:

(a) \ $|t_1|\le \gamma ^B<|t_0|, \gamma \in b^{s\restriction\eta }\in
t^*_1\longrightarrow t_0(\gamma ^B)=0$ or $s(\eta ).$

(b) \ $|t_1|\le \gamma ^C<|t_0|, \gamma =\lan\gamma _0,\gamma _1\ran,
\gamma _0\in A\longrightarrow t_0(\gamma ^C)=0.$

An $R^s$-generic is determined by a function $T:\ \alpha ^+\longrightarrow 2$
such that $s(\eta )=0$ iff $T(\gamma ^B)=0$ for sufficiently large $\gamma
\in b^{s\restriction \eta } $ and such that for $\gamma _0<\alpha ^+:$
$\gamma _0\in A$ iff $T(\gamma ^C)=0$ for sufficiently large $\gamma =\lan\gamma _0,\gamma _1\ran<\alpha ^+.$

Now we come to the definition of the limit coding, which incorporates the
idea of ``coding delays''. Suppose $s\in S_\alpha ,$ $\alpha \in\Card'$ and
$\vec p=\lan(p_\beta ,p^*_\beta )|\beta \in \Card \cap \alpha \ran$ where
$p_\beta \in S_\beta $ for each $\beta \in\Card\cap\alpha.$ We wish to define:
``$\vec p$ codes $s$''. A natural definition would be: \ for
$\eta <|s|,$ $p_\beta (f^{s\restriction\eta }(\beta ))=s(\eta )$ for
sufficiently large $\beta \in \Card \cap \alpha.$  There are problems with
this definition however. First, to avoid conflict with the successor coding
we should use $f^{s\restriction\eta }(\beta )^D$ instead of
$f^{s\restriction\eta }(\beta ).$  And it is convenient and sufficient to
only require the above for $\beta \in \Card^+\cap\alpha.$  However,
there are still serious difficulties in making sure that the coding of $s$
is consistent with the codings of $p_\beta $ by $\vec p\restriction \beta
,$ for $\beta \in \Card'\cap\alpha.$  To solve these problems Jensen used
$\square$ to make these codings almost disjoint, for singular $\alpha;$
this creates new difficulties, resulting from the fact that the singular
and inacessible codings are thereby different.

We introduce Coding Delays to facilitate an easier proof of extendibility
of conditions. The rough idea is to code $s(\eta )$ not at $f^{s\restriction
\eta }(\beta )^D$ but instead just after the least ordinal $\ge
f^{s\restriction \eta }(\beta )^D$ where $p_\beta $ takes the value $1.$

\vskip10pt

\noindent
{\bf Definition.}\ Suppose $\alpha \in \Card',$ $s\in S_\alpha.$
Let $\tilde \mu ^s$ be defined just like $\mu ^s$ but with the
requirement ``$\mu '\mu =\mu $ for $\mu '<\mu $'' replaced by the
weaker requirement ``$\mu $ a limit ordinal.'' Then note that
$\widetilde{\cal{A}}^s=L_{\tilde\mu^s}[A\cap \alpha ,s]$ belongs to $\cal{A},$
contains $s$ and the $\Sigma _1$ Hull $(\alpha \cup \{A\cap\alpha ,s\})$ in
$\widetilde{\cal{A}}^s=\widetilde{\cal{A}}^s.$ Now $X$ {\it codes}  $s$
if $X$ is the $\Sigma _1$ theory of $\tilde{\cal{A}}^s$ with
parameters from $\alpha \cup \{A\cap\alpha ,s\}$ (viewed as a subset of
$\alpha$).

\vskip10pt

\noindent
{\bf Definition.}\  (Limit Coding) \ Suppose $s\in S_\alpha ,\alpha \in
\Card'$  and $\vec p=\lan(p_\beta ,p^*_\beta )|\beta \in \Card\cap\alpha \ran$
where     $p_\beta \in S_\beta $ for each $\beta \in \Card\cap\alpha.$ We wish
to define ``$\vec p$ codes $s$''.  First we define a sequence
$\lan s_\gamma |\gamma \le\gamma _0\ran$ of elements of $S_\alpha $ as follows.
Let $s_0=\phi _\alpha.$  For limit $\gamma \le\gamma _0,$ $s_\gamma
=\cup\{s_\delta |\delta <\gamma \}.$ Now suppose $s_\gamma $ is defined and  let
$f^{s_{\gamma }}_p(\beta )=$ least $\delta \ge f^{s_{\gamma }}(\beta )$
such that $p_\beta (\delta ^D)=1,$ if such a $\delta $ exists. If
$f^{s_{\gamma }}_{\vec{p}}(\beta )$ is undefined for cofinally many $\beta
\in\Card^+\cap\alpha $ then set $\gamma _0=\gamma.$  Otherwise define
$X\subseteq \alpha $ by: \ $\delta \in X$ iff $p_\beta \bigl((f^{s_{\gamma
}}_{\vec p}(\beta )+1+\delta )^D\bigr)=1$ for sufficiently large $\beta \in
\Card^+\cap\alpha.$ If  Even $(X)$ codes an element $t$ of $S_\alpha $
extending $s_\gamma $ such that $f^{s_\gamma }_{\vec{p}},X\in\cal{A}^t$
then set $s_{\gamma +1}=t.$ Otherwise let $s_{\gamma +1}$ be $s_\gamma
*X^E$ if this definition yields $f^{s_{\gamma
}}_{\vec{p}}\in\cal{A}^{s_{\gamma +1}}$ (and otherwise $\gamma _0=\gamma ).$
Now $\vec{p}$ {\it exactly codes} $s$ if $s=s_\gamma $ for some $\gamma
\le\gamma _0$ and $\vec{p}$ {\it codes} $s$ if $s\le s_\gamma $ for some
$\gamma \le \gamma _0.$

\vskip10pt

\noindent
{\bf Definition} \ (The Conditions) \ A {\it condition} in $P$ is a
sequence $p=\lan(p_\alpha ,p^*_\alpha )|\alpha \in \Card,\alpha \le\alpha (p))$
where $\alpha (p)\in\Card$ and:

(a) \ $p_{\alpha (p)}\in S_{\alpha (p)}, p^*_{\alpha (p)}=\phi .$

(b) \ For $\alpha \in\Card \cap\alpha (p), $ $(p_\alpha ,p^*_\alpha )\in
R^{p_{\alpha ^+}}.$

(c) \ For $\alpha \in\Card',$ $\alpha \le\alpha (p),$ $p\restriction\alpha
\in\cal{A}^{p_{\alpha }},$ $p\restriction\alpha $ exactly codes $p_\alpha .$

(d) \ For $\alpha \in \Card',$ $\alpha \le\alpha (p),$ $\alpha $ inaccessible
in $\cal{A}^{p_\alpha},$ there exists CUB $C\subseteq \alpha ,$
$C\in\cal{A}^{p_\alpha }$ such that $\beta \in C\longrightarrow p^*_\beta
=\phi .$

Conditions are ordered by: \ $p\le q$ iff $\alpha (p)\ge\alpha (q),$
$p(\alpha )\le q(\alpha )$ in $R^{p_{\alpha ^+}}$ for $\alpha
\in\Card\cap\alpha (p)\cap(\alpha (q)+1)$ and $p_{\alpha (p)}$ extends
$q_{\alpha (p)}$ if $\alpha (q)=\alpha (p).$

It is also useful to define some approximations to $P:$ \ For $\alpha \in
\Card, P^{<\alpha }$ denotes the set of all conditions $p$ such that
$\alpha (p) <\alpha .$ Also for $s\in S_\alpha ,$ $\omega <\alpha \in\Card,$
$P^s$ denotes $P^{<\alpha }$ together with all $p\restriction\alpha $ for
conditions $p$ such that $\alpha (p)=\alpha ,$ $p_{\alpha (p)}\le s.$
To order conditions in $P^s$, first define $p^+=p$ for $p\in P^{<\alpha }$
and for $p\in P^s-P^{<\alpha},$ $p^+\restriction \alpha =p$ and $p^+(\alpha
)=(s\restriction\eta ,\phi ),$ $\eta $ least such that $p\in
P^{s\restriction\eta };$ then $p\le q$ iff $p^+\le q^+$ as conditions in
$P.$

It is worth noting that (c) above implies that $f^{p_\alpha }$ dominates
the coding of $p_\alpha $ by $p\restriction\alpha ,$ in the sense that
$f^{p_\alpha }$ strictly dominates each $f^{p_\alpha \restriction\eta
}_{p\restriction\alpha}$, $\eta <|p_\alpha |$ on a tail of
$\Card^+\cap\alpha .$ The purpose of (d) is to guarantee that extendibility
of conditions at (local) inaccessibles is not hindered by the Successor
Coding (see the proof of Extendibility below).

We now embark on a series of lemmas which together show that $P$ is the
desired forcing: $P$ preserves cofinalities and if $G$ is $P$-generic over
$\lan L[A], A\ran$ then $L[A,G]=L[X]$ for some $X\subseteq \omega _1,$ $A$
is $L[X]$-definable from the parameter $X.$

\vskip10pt

\noindent
{\bf Lemma 5} \ (Distributivity for $R^s$) \ Suppose $\alpha \in \Card,
s\in S_{\alpha ^+}.$ Then $R^s$ is $\alpha ^+$-distributive in $\cal{A}^s$:
if $\lan D_i|i<\alpha \ran\in \cal{A}^s$ is a sequence of dense subsets
of $R^s$ and $p\in R^s$ then there is $q\le p$ such that $q$ meets each
$D_i.$

\vskip10pt

\noindent
{\bf Proof} \ Choose $\mu <\mu ^s$ to be a large enough limit ordinal such
that $p,\lan D_i|i<\alpha \ran$, $\cal{A}^{<s}\in\cal{A}=L_\mu [A\cap
\alpha ^+,s].$ Let $\lan\alpha _i|i<\alpha \ran$ enumerate the first
$\alpha $ elements of $\{\beta <\alpha ^+|\beta =\alpha ^+\cap \Sigma _1$
Hull of $(\beta \cup\{p,\lan D_i|i<\alpha \ran, \cal{A}^{<s}\})$ in
$\cal{A}\}.$

Now write $p$ as $(t_0,t^*_0)$ and successively extend to $(t_i,t^*_i)$ for
$i\le\alpha $ as follows: $(t_{i+1}, t_{i+1})$ is the least extension of
$(t_1,t^*_i)$ meeting $D_i$ such that $t_{i+1}^*$ contains
$\{b^{s\restriction \eta }|\eta \in H_i\cap|s|\}$ where $H_i=\Sigma _1$ Hull
of $\alpha _i\cup\{p,\lan D_i|i<\alpha \ran,\cal{A}^{<s}\}$ in $\cal{A}$
and: \ (a) \ If $b^{s\restriction\eta } \in t_i^*,$ $s(\eta )=1$ then
$t_{i+1}(\gamma ^\beta )=1$ for some $\gamma \in b^{s\restriction\eta },$
$\gamma >|t_i|.$ \ (b) \ If $\gamma _0\notin A,$ $\gamma _0<|t_i|$ then
$t_{i+1}(\lan\gamma _0,\gamma _1\ran^C)=1$ for some $\gamma _1>|t_i|.$

The lemma reduces to:

\vskip10pt

\noindent
{\bf Claim} \ $(t_{\lambda },t^*_\lambda )=$ greatest lower bound to $\lan
(t_i,t_i^*)|i<\lambda \ran$ exists for limit $\lambda \le \alpha .$

\vskip10pt

\noindent
{\bf Proof of Claim.} \ We must show that $t_\lambda =\cup\{t_i|i<\lambda \}$
belongs to $S_\alpha.$  Note that $\lan t_i|i<\lambda \ran$ is definable
over $\overline{H}_\lambda =$ transitive collapse of $H_\lambda $ and by
construction, $t_\lambda $ codes $\overline{H}_\lambda $ definably over
$L_{\bar\mu _\lambda }[t_\lambda ],$ where $\bar\mu _\lambda =$ height of
$\overline{H}_\lambda .$ So $t_\lambda $ is reshaped, as $|t_\lambda |$ is
singular, definably over $L_{\bar\mu _\lambda }[t_\lambda ].$ By Lemma 1,
$\bar\mu _\lambda <(|t_\lambda |^+)^L$ if $\alpha \ge\omega _2.$ So
$t_\lambda $ belongs to $S_\alpha.$ $\dashv$

The next lemma illustrates the use of coding delays:

\vskip10pt

\noindent
{\bf Lemma 6} \ (Extendibility for $P^s$) \ Suppose $p\in P^s,$ $s\in
S_\alpha, $ $X\subseteq \alpha ,$ $X\in\cal{A}^s.$ Then there exists $q\le
p$ such that $X\cap \beta \in \cal{A}^{q_\beta }$ for each $\beta \in
\Card\cap\alpha.$

\vskip10pt

\noindent
{\bf Proof} \ Let $Y\subseteq \alpha $ be chosen so that Even $(Y)$ codes $s$
and Odd $(Y)$ is the $\Sigma _1$ theory of $\cal{A}$ with parameters from
$\alpha \cup \{A\cap \alpha ,s\},$ where $\cal{A}$ is an initial segment of
$\cal{A}^s$ large enough to extend $\cal{\widetilde A}^s$ and to contain $X,p$.
For $\beta \in\Card\cap\alpha ,$ let $\cal{\overline A}_\beta=$ transitive
collapse of 
$\Sigma _1$ Hull $(\beta \cup\{A\cap \alpha ,s\})$ in $\cal{A},$ and $g(\beta) = {\beta}^+$ of $\cal{\overline A}_\beta.$ 

Define $q$ as follows: $q_\beta =s_\beta $ if Even $(Y\cap\beta )$ codes
$s_\beta \in S_\beta ,$ $q_\beta =p_\beta *(Y\cap\beta )^E$ for other
$\beta \in \Card'\cap\alpha ,$ $q_\beta =p_\beta *\vec O*1*(Y\cap \beta
)^D$ where $\vec{O}$ has length $g(\beta )$ for $\beta \in
\Card^+\cap\alpha .$ And $q^*_\beta =p^*_\beta $ for all $\beta \in \Card
\cap\alpha.$

As $g\restriction\beta ,$ $Y\cap\beta $ are definable over $\cal{\overline
A}_\beta $ for $\beta \in\Card\cap\alpha $ we get $g\restriction\beta ,$
$Y\cap\beta \in\cal{A}^{s_\beta }$ when Even$(Y\cap\beta )$ codes $s_\beta
\in S_\beta .$ Also $g\restriction\beta ,$ $Y\cap\beta \in
\cal{A}^{q_{\beta }}$ for other $\beta \in \Card'\cap\alpha $ as Odd
$(Y\cap\beta )$ codes $\cal{\overline A}_\beta.$ And note that for all
$\beta \in \Card'\cap\alpha ,$ $g\restriction \beta $ dominates
$f^{p_{\beta }}$ on a final segment of $\Card^+\cap\beta ,$ unless
Even $(Y\cap\beta )$ codes $s_\beta =p_\beta ,$ in which case
$q\restriction\beta $ exactly codes $s_\beta $ because $p\restriction\beta
$ does.

So we conclude that $q\restriction\beta $ exactly codes $q_\beta $ for
sufficiently large $\beta \in \Card'\cap\alpha $ and clearly $X\cap \beta
\in \cal{A}^{q_\beta }$ for such $\beta.$ Apply induction on $\alpha $ to
obtain this for all $\beta \in\Card'\cap\alpha.$ Finally, note that the
only problem in verifying $q\le p$ is that the restraint $p^*_\beta $ may
prevent us from making the extension $q_\beta $ of $p_\beta $ when $q_\beta
=s_\beta ,$ Even $(Y\cap \beta )$ codes $s_\beta.$ But property (d) in the
definition of condition guarantees that $p^*_\beta =\phi$ for $\beta $ in a
CUB $C\subseteq \alpha ,$ $C\in \cal{A}^{s}.$ We may assume that $C\in
\cal{A}$ and hence for sufficiently large $\beta $ as above we get $\beta \in
C$ and hence $p^*_\beta =\phi.$ So $q\le p$ on a final segment of
$\Card\cap\alpha ,$ and we may again apply induction to get $q\le p$
everywhere. $\dashv$

\vskip10pt

The key idea of Jensen's proof lies in the verification of distributivity
for $P^s$. Before we can state and prove this property we need some definitions.

\vskip10pt

\noindent
{\bf Definition} \ Suppose $\beta \in \Card^+\cap\alpha $ and $D\subseteq
P^s,$ $s\in S_\alpha.$ $D$ is $\beta$-{\it dense} on $P^s$ if $\forall p\in
P^s\exists q\in P^s(q\le p, q$ meets $D$ and $q\restriction \beta
=p\restriction \beta).$  $X\subseteq \Card \cap\alpha$ is {\it thin in}
$\cal{A}^s$ if $X\in \cal{A}^s$ and for each inaccessible $\beta \le \alpha,$
$\cal{A}^s\vDash X\cap \beta$ is not stationary in $\beta.$ A function
$f: \ \Card\cap\alpha \longrightarrow V$ in $\cal{A}^s$ is {\it small in}
$\cal{A}^s$ if for each $\beta \in\Card\cap\alpha ,$ $f(\beta )\in
H^{\cal{A}^s} _{\beta ^{++}},$ $\card(f(\beta ))\le \beta $ in $\cal{A}$
and Support $(f)=\{\beta \in\Card\cap\alpha |f(\beta )\neq\phi \}$ is
thin in $\cal{A}^s$. If $D\subseteq P^s$ is predense and $p\in P^s,$ $\beta
\in\Card$ we say that $p$ {\it reduces D below} $\beta $ if for some
$\gamma \in \Card^+$ $\gamma \le\beta,$ $\{r|r\cup p\restriction
[\gamma ,\alpha)$ meets $D\}$ is predense on $P^{p_\gamma }$ below
$p\restriction \gamma.$  Finally, for $p\in P^s,f$ small in $\cal{A}^s$
we define $\Sigma ^p_f=$ all $q\le p$ in $P^s$ such that whenever $\beta
\in \Card\cap\alpha,$ $D\in f(\beta ),D$ predense on $P^{p_{\beta ^+}}$
then $q$ reduces $D$ below $\beta.$

\vskip10pt

\noindent
{\bf Lemma 7} \  (Distributivity for $P^s$) \ Suppose $s\in S_{\beta ^+},$ $\beta \in\Card.$

(a) \ If $\lan D_i|i<\beta \ran\in\cal{A}^s; D_i$ $i^+$-dense on $P^s$ for
each $i<\beta $ and $p\in P^s$ then there is $q\le p,$ $q$ meets each $D_i.$

(b) \ If $p\in P^s,f$  small in $\cal{A}^s$ then there exists $q\le p,$
$q\in \Sigma ^p_f.$

\vskip10pt

\noindent
{\bf Proof} \ We demonstrate (a) and (b) by a simultaneous induction on
$\beta.$ If $\beta =\omega $ or belongs to $\Card^+$ then by induction (a)
reduces to the $\beta ^+$-distributivity of $R^s$ in $\cal{A}^s,$ Lemma 5.
And (b) reduces to: \ if $S$ is a collection of $\beta $-many predense
subsets of $P^s,$ $S\in\cal{A}^s$ then $\{q\in P^s|q$ reduces each $D\in S$
below $\beta \}$ is dense on $P^s.$ Again this follows from Lemma 5 since
$P^s$ factors as $R^s*Q$ where $1^{R^s}\Vdash Q$ is $\beta ^+-cc,$
and hence any $p\in P^s$ can be extended to $q\in P^s$ such that
$D^q=\{r\in D|q(\beta )\le r(\beta )$ in $R^s\}$ is predense $\le q$ for each
$D\in S$ and hence $q$ reduces each $D\in S$ below $\beta.$

Now suppose that $\beta $ is inaccessible. We first show that (b) holds for
$f,$ provided $f(\beta )=\phi.$ First select a CUB $C\subseteq \beta $ in
$\cal{A}^s$ such that $\gamma \in C\longrightarrow f(\gamma )=\phi $ and
extend $p$ so that $f\restriction \gamma , C\cap \gamma $ belong to
$\cal{A}^{p_\gamma }$ for each $\gamma \in \Card\cap \beta ^+.$ Then we can
successively extend $p$ on $[\beta _i^+,\beta _{i+1}]$ in the least way so as 
to meet $\Sigma ^p_f$ on $[\beta _i^+, \beta _{i+1}]$, where $\lan\beta
_i|i<\beta \ran$ is the increasing enumeration of C. At limit stages
$\lambda ,$ we still have a condition, as the sequence of first $\lambda $
extensions belongs to $\cal{A}^{p_{\beta_{\lambda }}}.$ The final
condition, after $\beta $ steps, is an extension of $p$ in $\Sigma ^p_f.$

Now we prove (a) in this case. Suppose $p\in P^s$ and $\lan D_i|i<\beta
\ran\in\cal{A}^s,$ $D_i$ is $i^+$-dense on $P^s$ for each $i<\beta.$
Let $\mu _0<\mu ^s$ be a big enough limit ordinal so that
$\lan D_i|i<\beta \ran,p,$ $\tilde\mu ^s\in L_{\mu _0}$ $[A\cap \beta
^+,s]$ and for $i<\beta $ let $\mu _i=\mu _0+\omega \cdot i<\mu ^s.$ For
any $X$ we let $H_i(X)$ denote $\Sigma _1$ Hull$(X\cup\{\lan D_i|i<\beta \ran,$
$p,\tilde\mu ^s,s,$ $A\cap\beta ^+\})$ in $L_{\mu _i}[A\cap \beta ^+,s].$

Let $f_i: \ \Card\cap \beta \longrightarrow V$ be defined by:
$f_i(\gamma )=H_{\gamma ^{++}}\cap H_i(\gamma )$ if $i<\gamma \in
H_i(\gamma ),$ $i<\gamma <\beta $ and $f_i(\gamma )=\phi $ otherwise. Then
each $f_i$ is small in $\cal{A}^s$ and we inductively define $p=p^0\ge
p^1\ge\dots$ in $P^s$ as follows: $p^{i+1}=$ least $q\le p^i$ such that:

(a) \ $q(\beta )$ meets all predense $D\subseteq R^s,$ $D\in H_i(\beta ).$

(b) \ $q$ meets $\Sigma ^{p^i}_{f_i}$ and $D_i.$

(c) \ $q\restriction i^+=p^i\restriction i^+.$

For limit $\lambda \le\beta $ we take $p^\lambda $ to be the greatest
lower bound to $\lan p^i|i<\lambda \ran,$ if it exists.

\vskip10pt

\noindent
{\bf Claim} \ $p^\lambda $ is a condition in $P^s,$ where $p^\lambda
(\gamma )=(\cup\{p^i_\gamma |i<\lambda \},$ $\cup\{p^i_\gamma{}^*|i<\lambda
\})$ for each $\gamma \in\Card\cap\beta ^+.$

First we verify that $p^\lambda _\gamma =\cup\{p^i_\gamma
|i<\lambda \}$ belongs
to $S_\gamma.$ Let $\overline{H}_\lambda (\gamma )$ be the transitive
collapse of $H_\lambda (\gamma )$ and write $\overline{H}_\lambda (\gamma
)$ as $L_{\bar\mu }[\overline{A},\bar{s}],$ $\overline{P}=$ image of
$P^s\cap H_\lambda (\gamma )$ under transitive collapse, $\bar\beta =$ image
of $\beta $ under collapse. Also write $\overline{P}$ as $R^{\bar s}*
P^{\bar G_{\bar\beta }}$ where $\overline{G}$ denotes an $R^{\bar
s}$-generic (just as $P^s$ factors as $R^s*P^{G_{\beta }},$ $G_\beta $
denoting an $R^s$-generic).

Now the construction of the $p^i$'s (see conditions (a), (b)) was designed
to guarantee that if $\gamma \in H_\lambda (\gamma )$ then
$\overline{G}_{\bar\beta }=\{\bar p\in R^{\bar s}|\bar p$ is extended by some
$\bar p^i(\bar\beta )\}$ is $R^{\bar s}$-generic over $\overline{H}_\lambda
(\gamma ),$ where $\bar p^i$ = image of $p^i$ under collapse, and that for
each $\gamma <\bar\delta <\bar\beta$ in $\Card^+(\overline{H}_\lambda
(\gamma )),$ $\{\bar p|\bar p$ is extended by some $\bar p^i\restriction
[\gamma ,\bar\delta )$ in $\overline{P}_\gamma ^{\bar p^{i}_{\bar\delta }}\}$
is $\overline{P}_\gamma ^{\overline{G}_{\bar\delta }}$-generic over
 $\cal{A}^{<\bar G_{\bar\delta }}=\cup\{\cal{A}^{<\bar
p^{i}_{\bar\delta }} |i<\lambda \}$ where
$\overline{P}_\gamma ^{\bar p^{i}_{\bar\delta }}$ denotes the image under
collapse of $P^{p^{i}_{\delta }}_\gamma =\{q\restriction[\gamma ,\delta
)|q\in P^{p^{i} _{\delta }}\}, \bar\delta =$ image of $\delta $ under collapse.

\vskip10pt

\noindent
{\bf Note:} \ We do {\it not} necessarily have the previous claim for
$\bar\delta =\bar\beta ,$ and this is the source of our need for $\sim
O^\#$ in this proof.

By induction, we have the distributivity of $P^t$ for $t\in S_\delta ,$
$\delta \in \Card^+\cap\beta ,$ and hence that of $\overline P^{\bar t}$ for
$\bar t\in \overline{S}_{\bar\delta }, \bar\delta
\in\Card^+(\overline{H}_\lambda (\gamma )),$ $\bar\delta <\bar\beta.$
So the ``weak'' genericity of the preceding paragraph implies that:

(d) \ $L_{\bar\beta }[A\cap\gamma ,p^\lambda _\gamma ]\vDash p^\lambda
_\gamma |$ is a cardinal.

\noindent
Also:

(e) \ $L_{\bar\mu }[A\cap\gamma ,p^\lambda _\gamma ]\vDash|p^\lambda
_\gamma |$ is $\Sigma _1$-singular.

\noindent
Thus $p^\lambda _\gamma \in S_\gamma$ (by (e)) provided we can show that when
$\gamma \ge\omega _2,$ $\bar\mu <(|p^\lambda _\gamma |^+)^L.$ But
$\overline{H}_\lambda (\gamma )\stackrel{\sim}{\longrightarrow} H_\lambda
(\gamma )$ gives a $\Sigma _1$-elementary embedding with critical point
$|p^\lambda _\gamma|,$ so by Lemma 1, this is true. Also note that we now get
$p^\lambda \restriction \gamma \in\cal{A}^{p^{\lambda }_{\gamma }}$ as well,
since $p^\lambda \restriction\gamma $ is definable over
$\overline{H}_\lambda (\gamma )$ and we defined $\cal{A}^{p^{\lambda
}_{\gamma }}$ to be large enough to contain $\overline{H}_\lambda (\gamma
),$ since $L_{\bar\beta }\vDash |p^\lambda _\gamma |$ is a cardinal by (d).

The previous argument applies also if $\gamma =\beta,$ using the
distributivity of $R^s,$ or if $\gamma =\beta \cap H_\lambda (\gamma )$, using
the fact that $p^\lambda _\beta $ collapses to $p^\lambda _\gamma.$ If
$\gamma <\gamma ^*=\operatorname{min} (H_\lambda (\gamma )\cap[\gamma
,\beta ))$ then we can apply the first argument to get the result for
$\gamma ^*,$  and then the second argument to get the result for $\gamma.$

Finally, to prove the Claim we must verify the restraint condition (d) in
the definition of $P.$ Suppose $\gamma $ is inaccessible and for $i<\lambda
$ let $C^i$ be the least CUB subset of $\gamma $ in $\cal{A}^{p^i_\gamma }$
disjoint from $\{\bar\gamma <\gamma |p^i_{\bar\gamma}*\neq\phi \}.$ If
$\lambda <\gamma $ then $\bigcap \{C^i|i<\lambda \}$ witnesses the restraint
condition for $p^\lambda $ at $\gamma $, if $\gamma < \lambda$ then the 
restraint condition for $p^{\lambda}$ at $\gamma$ follows by induction on $\lambda $
and if $\gamma =\lambda $ then $\Delta \{C^i|i<\lambda \}$ witnesses the
restraint condition for $p^\lambda $ at $\gamma ,$ where $\Delta $ denotes
diagonal intersection.

Thus the Claim and therefore (a) is proved in case $\beta $ is
inaccessible. To verify (b) in this case, note that as we have already
proved (b) when $f(\beta )=\phi $ it suffices to show: if $\lan D_i|i<\beta
\ran\in\cal{A}^s$ is a sequence of dense subsets of $P^s$ then $\forall
p\exists q\le p$ ($q$ reduces each $D_i$ below $\beta).$  But using
distributivity we see that $D^*_i = \{q|q$ reduces $D_i$ below $i^+\}$ is
$i^+$-dense for each $i<\beta $ so again by distributivity there is $q\le p$
reducing each $D_i$ below $i^+.$

We are now left with the case where $\beta $ is singular. The proof of (a)
can be handled using the ideas from the inaccessible case, as follows.
Choose $\lan\beta _i|i<\lambda _0\ran$ to be a continuous and cofinal
sequence of cardinals $<\beta ,\lambda _0<\beta _0.$ First, we argue that
$p\in P^s$ can be extended to meet $\Sigma ^p_f$ for any $f$ small in
$\cal{A}^s$, provided $f(\beta )=\phi :$ Extend $p$ if necessary so that
for each $\gamma \in \Card\cap \beta ^+,$ $f\restriction\gamma $ and
$\{\beta _i|\beta _i<\gamma \}$ belong to $\cal{A}^{p_{\gamma }}$.  Now
perform a  construction like the one used to prove distributivity in the
inacessible case, extending $p$ successively on $[\beta _0,\beta _i^+]$ so
as to meet $\Sigma ^p_f$ on $[\beta _0,\beta ^+_i]$ as well as appropriate
$\Sigma ^{p_i}_{f_i}$'s defined on $[\beta _0,\beta _i^+]$ to guarantee
that $p^\lambda $ is a condition for limit $\lambda \le\lambda _0$. Note
that each extension is made on a bounded initial segment of $[\beta _0,\beta )$
and therefore by induction $\Sigma ^p_f,\Sigma ^{p_i}_{f_i}$ can be met on
these intervals. The result is that $p$ can be extended to meet $\Sigma
^p_f$ on a final segment of $\Card\cap\beta $ and therefore by induction
can be extended to meet $\Sigma ^p_f.$ Second, use the density of $\Sigma
^p_f$ when $f(\beta )=\phi $ to carry out the distributivity proof as we
did in the inaccessible case. And again, (b) follows from (a). This
completes the proof of Lemma 7. $\dashv$

\vskip10pt

Now the same argument as used above also shows:

\vskip10pt

\noindent
{\bf Lemma 8} \ (Distributivity for $P$) \ If $\lan D_i|i<\beta \ran$ is
$\lan L[A],A\ran$-definable, each $D_i$ is $i^+$-dense on $P$ and $p\in P$
then there exists $q\le p,$ $q$ meets each $D_i.$

Extendibility for $P^s$ and Distributivity for $P$ give us the conclusions
of Theorem 4. This completes the proof.

\vskip20pt

\begin{center} {\bf References} \end{center}

\vskip10pt

\noindent
Beller-Jensen-Welch [82] \ {\it Coding the Universe,} Cambridge University
Press.

\vskip5pt

\noindent
Friedman [94] \ A Simpler Proof of Jensen's Coding Theorem, Annals of Pure
and Applied Logic, vol. 70, No.1, pages 1--16.

\vskip5pt

\noindent
Jensen-Solovay [68] \ Some Applications of Almost Disjoint Sets, in {\it
Mathematical Logic and the Foundations of Set Theory,} North Holland,
pages 84--104.

\end{document}